\documentclass[12pt]{amsart}
\usepackage{euscript,oldgerm,geometry}
\usepackage{amssymb,amsmath}
\usepackage{color} 
 
\newcommand{\conju}[1]{\overline{#1}}
\newcommand{\refe}[1]{(\ref{#1})}
\newcommand{\dst}{\displaystyle}

\newcommand{\PP}{{\mathbb P}}
\newcommand{\CC}{{\mathbb C}}

\newcommand{\pe}[2]{\langle#1,#2\rangle}
\renewcommand{\u}{\mbox{${\textswab u}$}}
\newcommand{\newu}{\mbox{$\widetilde{\textswab u}$}}
\newcommand{\p}{\widetilde{P}}

\newcommand{\K}{\mathrm{K}}

\newcommand{\vp}{\widetilde{{P}}_n}

\newcommand{\half}{\mbox{$\frac12$}}

\newcommand{\bq}{\begin{equation}}
\newcommand{\eq}{\end{equation}}
\newcommand{\ba}{\begin{array}}
\newcommand{\ea}{\end{array}}

\newtheorem{theorem}{Theorem}
\newtheorem{rem}[theorem]{Remark}

\begin{document}

\title{On the generalized Askey-Wilson Polynomials}

\author{R. \'{A}lvarez-Nodarse}
\address{IMUS \& Departamento de An\'alisis Matem\'atico, Universidad de Sevilla.
Apdo. 1160, E-41080  Sevilla, Spain}\email{ran@us.es}

\author{R. Sevinik Ad\i g{\"{u}}zel}
\address{Department of Mathematics, Faculty of Science, Sel\c{c}uk University, 42075,
Konya, Turkey}\email{sevinikrezan@gmail.com}



\begin{abstract}
In this paper a generalization of Askey-Wilson polynomials is introduced.
These polynomials are obtained from the Askey-Wilson polynomials via the addition of two mass points
to the weight function  of them at the points $\pm1$.
Several properties of such new family are considered, in particular the three-term recurrence 
relation and the representation as basic hypergeometric series.
\end{abstract}

\maketitle 

\thanks{\centerline{\small \textit{Dedicated  to  Paco Marcell\'an on  the  occasion  of his  60th  birthday}}}

\section{Introduction}

The Krall-type polynomials are orthogonal with respect to a
linear functional $\newu$ obtained from a
quasi-definite functional $\u:\mathbb{P}\mapsto\CC$  ($\mathbb{P}$, denotes the space
of complex polynomials with complex coefficients) via the addition of delta Dirac measures.
These polynomials appear as eigenfunctions of a fourth order 
linear differential operator with polynomial coefficients
that do not depend on the degree of the polynomials. 
They were firstly considered by Krall in \cite{hkrall} 
(for a more recent reviews see \cite{a-nmp} and \cite[chapter XV]{akrall}). 
In fact, H.\ L.\ Krall discovered that there  are only three extra families of orthogonal 
polynomials apart from the classical polynomials of Hermite, Laguerre and Jacobi
that satisfy such a fourth order differential equation which are
orthogonal with respect to measures that are not absolutely continuous with
respect to the Lebesgue measure. This result motivated the study of the polynomials orthogonal 
with respect to the more general weight functions \cite{koko, koo84} that could contain
more instances of orthogonal polynomials being eigenfunctions of higher-order 
differential equations \cite[chapters XVI, XVII]{akrall}.

In the last years the study of such polynomials have been considered by many authors
(see e.g. \cite{RJF,RJ,HG96,hai01,kwo00,VZ01} and the references therein)
with a special emphasis on the case when the starting functional $\u$ is a classical
continuous, discrete or $q$-linear functional with the linear type lattices 
(for more details see  \cite{RRC, RJ} and references therein). In fact, for the 
$q$-case some examples related with the $q$-Laguerre and the little $q$-Jacobi polynomials
were constructed by Haine and Gr\"unbaum in \cite{HG96} using the Darboux transformation. 
Later on, in \cite{VZ01}, Vinet and Zhedanov presented a more complete study for the little 
$q$-Jacobi polynomials. In these both cases, the $q$-Krall polynomials satisfy a 
higher order $q$-difference equations with polynomial coefficients independent of $n$. 
For the discrete case the problem was solved very recently by A. Dur\'an using 
a new method (see \cite{dur1,dur2}, for details).

For the general $q$-quadratic lattice only few results were known. 
An important contribution to this case was done in \cite{HP} where
the authors considered a generalized Askey-Wilson polynomials 
by adding mass points. They showed that the resulting orthogonal polynomials satisfy 
a higher order $q$-difference equation with polynomial coefficients independent of $n$,
only if the masses are added at very specific points out of the interval of orthogonality 
$[-1,1]$. Another contribution to this problem was done in \cite{RR}, where
a general theory of the the Krall-type polynomials on non-uniform lattices was developed. 
In fact, in \cite{RR} the authors studied the polynomials $\p_n(s)_q$ which are orthogonal with respect
to the linear functionals $\newu=\u+\sum_{k=1}^N A_k\delta_{x_k}$ defined on the 
$q$-quadratic lattice $x(s)=c_1q^s+c_2q^{-s}+c_3$ and considered, as a representative 
example, the Krall-type Racah polynomials (see also \cite{RRa}). In fact, in
\cite[\S5]{RR}, we posed the problem of obtaining a generalization of the Askey-Wilson polynomials
by adding two mass points at the end of the interval of orthogonality, motivated by the results in \cite{HP}.

Thus our main aim here is to study the orthogonal polynomials obtained via the 
addition of two mass points at the end of the interval of orthogonality of the
Askey-Wilson polynomials. The structure of the paper is as follows. In Section 2, some
preliminary results on the Askey-Wilson polynomials are presented 
as well as the most general expression for the
kernels  on the $q$-quadratic lattice $x(s)=c_1q^s+c_2q^{-s}+c_3$.
Our main results are in section 3, where we introduce a detailed study of the
generalized Askey-Wilson polynomials obtained from the classical Askey-Wilson 
polynomials by adding two mass points at $\mp1$.


\section{Preliminary results}
Here we  include some results of the theory of orthogonal polynomials
on the non-uniform lattice (for more details see e.g., \cite{ran,NSU})
\bq\label{redgen}
x(s)=c_1q^s+c_2q^{-s}+c_3.
\eq
The polynomials on non-uniform lattices $P_n(s)_q:=P_n(x(s))$ 
are the polynomial
solutions of the second order linear difference equation (SODE) of hypergeometric type
\begin{equation}\label{rez8}
\begin{split}
&\qquad \qquad A_sy(s+1) +B_s  y(s)+C_sy(s-1)+\lambda_ny(s)=0, \\
& A_s=\dst\frac{\sigma(s)+\tau(s)\Delta x(s-\frac{1}{2})}{\Delta x(s)\Delta x(s-\frac{1}{2})},\,\,\,
C_s=\dst\frac{\sigma(s)}{\nabla x(s)\Delta x(s-\frac{1}{2})},\,\,\, B_s=-A_s-C_s,
\end{split}
\end{equation}
where $\sigma(s)$ and $\tau(s)$ are polynomials of degree at most $2$ and exactly $1$,
respectively, and $\lambda_n$ is a constant.
They are orthogonal with respect to
the linear functional $\u:\PP_q\mapsto\CC$, where $\PP_q$ denotes
the space of polynomials on the lattice \refe{redgen}, 

\bq
\label{def-fun}\langle\u, P_n P_m \rangle = \delta_{m\,n} d_n^2,\quad 
\langle\u,P\rangle = \int_{x_0}^{x_1} P(x)_q \rho(x) dx
\eq
where $\rho$ is the \textit{weight} function 
and $d_n^2:=\langle\u, P_n^2 \rangle $.

Since the polynomials $P_n(s)_q$ are orthogonal with respect to a linear functional, they
satisfy a three-term recurrence relation (TTRR) \cite{ran,cbook}
\begin{equation}\label{rez10}
x(s)P_n(s)_q=\alpha_nP_{n+1}(s)_q+\beta_nP_n(s)_q+\gamma_nP_{n-1}(s)_q, \quad n=0, 1, 2, ...,
\end{equation}
with the initial conditions $P_0(s)_q=1,  P_{-1}(s)_q=0$, and
also the differentiation formulas \cite[Eqs. (5.65) and (5.67)]{ran} (or
\cite[Eqs. (24) and (25)]{RYR}
\begin{equation}\label{rez11}
\sigma(s)\frac{\nabla P_n(s)_q}{\nabla x(s)}=\overline{\alpha}_n P_{n+1}(s)_q
+\overline{\beta}_n(s)P_n(s)_q,
\end{equation}
\begin{equation}\label{rez12}
\Phi(s)\frac{\Delta P_n(s)_q}{\Delta x(s)}=\widehat{\alpha}_n P_{n+1}(s)_q
+\widehat{\beta}_n(s)P_n(s)_q,
\end{equation}
where $\Phi(s)=\sigma(s)+\tau(s)\Delta x(s-\frac{1}{2})$, and
$$
\overline{\alpha}_n = \widehat{\alpha}_n= -\frac{\alpha_n \lambda_{2n}}{[2n]_q},
\quad \overline{\beta}_n(s)=\frac{\lambda_{n}}{[n]_q}\frac{\tau_n(s)}{\tau'_n},\quad
\widehat{\beta}_n(s)=\overline{\beta}_n(s)-\lambda_n \Delta x(s-\half).
$$
Notice that from \refe{rez12} and the TTRR \refe{rez10} the following useful relation follows
\bq
\label{rez12a}
P_{n-1}(s)_q=\Theta(s,n)P_n(s)_q+\Xi(s,n) P_n(s+1)_q,
\eq
where
$$
\Theta(s,n)=\frac{\alpha_n}{\widehat{\alpha}_n\gamma_n}\left[
\frac{\Phi(s)}{\Delta x(s)}\!-\!\frac{\lambda_{2n}}{[2n]_q}(x(s)\!-\!\beta_n)+
\widehat{\beta}_n(s)\right],\,\,
\Xi(s,n)=-\frac{\alpha_n}{\widehat{\alpha}_n\gamma_n}\frac{\Phi(s)}{\Delta x(s)}.
$$
Using the Christoffel-Darboux formula for the $n$-th
reproducing kernels 
\begin{equation*}
\begin{array}{l}
\K_n(x(s_1),x(s_2)):=\dst\sum_{k=0}^{n}{\dst\frac{P_k(s_1)_qP_k(s_2)_q}{d_k^2}}=\dst\frac{\alpha_n}{d_n^2}
\frac{P_{n+1}(s_1)_q P_n(s_2)_q-P_{n+1}(s_2)_qP_{n}(s_1)_q}
{x(s_1)-x(s_2)},
\end{array}
\end{equation*}
and the relations \refe{rez11} and \refe{rez12}, respectively, to eliminate  $P_{n+1}$, we 
obtain the following two expressions  
\begin{equation}\label{rez14}
\begin{split}
\K_{n}(x(s),x(s_0)) & =\dst\frac{\alpha_nP_n(s_0)_q}{\overline{\alpha}_nd_n^2}
 \left\{\frac{\overline{\beta}_n(s_0)-\overline{\beta}_n(s)}{x(s)-x(s_0)}
 P_n(s)_q \dst+\frac{\sigma(s)}{x(s)-x(s_0)}\frac{\nabla P_n(s)_q}{\nabla x(s)}\right\}\\
 &-\frac{\alpha_n}{\overline{\alpha}_nd_n^2}\frac{\sigma(s_0)}{x(s)-x(s_0)}
\frac{\nabla P_n(s_0)_q}{\nabla x(s_0)}P_n(s)_q,
  \end{split}
 \end{equation}

 \begin{equation}\label{rez15}
\begin{split}
\K_{n}(x(s), x(s_0))& =\dst\frac{\alpha_nP_n(s_0)_q}{\widehat{\alpha}_nd_n^2}
 \left\{\frac{\widehat{\beta}_n(s_0)-\widehat{\beta}_n(s)}{x(s)-x(s_0)}
 P_n(s)_q \dst+\frac{\Phi(s)}{x(s)-x(s_0)}\frac{\Delta P_n(s)_q}{\Delta x(s)}\right\}\\
 &-\frac{\alpha_n}{\widehat{\alpha}_nd_n^2}\frac{\Phi(s_0)}{x(s)-x(s_0)}
 \frac{\Delta P_n(s_0)_q}{\Delta x(s_0)}P_n(s)_q.
  \end{split}
\end{equation}
Let us mention here that the above two formulas generalize to an arbitrary value $s_0$ 
the Eqs. (9) and (10) obtained in \cite[page 184]{RR}.

Next, we introduce the Askey-Wilson polynomials defined 
by the following basic series \cite{AW} (for the definition and properties of
basic series see e.g. \cite{GR}) 
\bq
\begin{split}\dst
P_n(x(s))_q:= P_n(x(s),a,b,c,d|q) &= \dst
\frac{(ab,ac,ad;q)_n}{(2a)^n(abcdq^{n-1};q)_n}\\
&\times{}_{4}\varphi_3 \left(\ba{c} q^{-n},abcdq^{n-1}, aq^{s},
aq^{-s}  \\ ab,ac,ad  \ea
\,\bigg|\, q \,,\, q \right).
\end{split}
\label{pol-wil-nu}
\eq
Notice that the Askey-Wilson polynomials are defined on 
the lattice $x(s)=\frac{q^s+q^{-s}}{2}$, $q^s=e^{i\theta}$ \cite{ata92}, which is a
particular case of \refe{redgen} when $c_1=c_2=1/2$ and $c_3=0$.
Their main characteristics (see Eqs. \refe{def-fun}, \refe{rez10}, \refe{rez11}, \refe{rez12}) 
are given in Table \ref{tabla1}.

\begin{table}[ht!]
\caption{Main data of the monic Askey-Wilson
polynomials  \cite{ks}\label{tabla1}}
\begin{center}
\begin{scriptsize}
{\renewcommand{\arraystretch}{.25}
\begin{tabular}{|@{}c@{}| | @{}c@{}|}\hline
 &  \\
$P_n(s)$ & ${P}_n(x(s),a,b,c,d|q)  \,,\quad x(s) =\frac{q^s+q^{-s}}{2}$,\,\, $q^s=e^{i\theta}$ , 
\quad $\Delta x(s)=\frac{q-1}{2}[q^s-q^{-s-1}]$ \\
 &  \\
\hline\hline
 &  \\
$\rho(s)$ &  \mbox{$\dst\frac{(q,ab,ac,ad,bc,bd,cd;q)_{\infty}h(x,1)h(x,-1)h(x,q^{1/2})h(x,-q^{1/2})}
{2\pi\sqrt{1-x^2}(abcd;q)_{\infty}h(x,a)h(x,b)h(x,c)h(x,d)}, \, h(x,\alpha)\!=\!\prod_{k=0}^{\infty}
\left(1-2\alpha xq^k+\alpha^2q^{2k}\right)$}  \\
 &  \\
 &  \mbox{
$x_0\!=\!-1$, $x_1\!=\!1$, $a,b,c,d\in\mathbb{R}$ or complex conjugate pairs if $a,b,c,d\in\mathbb{C}$ and
$\max(|a|,|b|,|c|,|d|)<1$} \\
 &  \\
\hline
 &  \\
$\sigma(s)$ & \mbox{
$ q^{-4s}(q^s-a)(q^s-b)(q^s-c)(q^s-d) $} \\
 &  \\
\hline
 &  \\
$\Phi(s)$ & \mbox{
$ q^{4s}(q^{-s}-a)(q^{-s}-b)(q^{-s}-c)(q^{-s}-d) $} \\
 &  \\
\hline
 &  \\
$\tau(s) $ &  \mbox{
\begin{tabular}{r}
$\frac{4}{q^{1/2}-q^{-1/2}}(ab+ac+ad+bc+bd+cd) x(s)-\frac{2}{(q^{1/2}-q^{-1/2})}
(a+b+c+d)$   
\end{tabular}}\\
 &  \\
\hline
 &  \\
$\tau_n(s) $ &  \mbox{
\begin{tabular}{r}
$ \frac{4q^n}{q^{1/2}-q^{-1/2}}(ab\!+\!ac\!+\!ad\!+\!bc\!+\!bd\!+\!cd)
x(s\!+\!\frac n2)-\frac{2q^{n/2}}{(q^{1/2}-q^{-1/2})}
(a+b+c+d)$
\end{tabular}}\\
 &  \\
 \hline
 &  \\
$\lambda _n$ &\mbox{
$4q^{-n+1}(1-q^n)(1-abcdq^{n-1})$} \\
 &  \\
\hline
 \hline
 &  \\
$d_n^2$ &    \mbox{$\dst\frac{2^{-2n}(q,ab,ac,ad,bc,bd,cd,abcdq^{2n};q)_{\infty}}
{(abcdq^{n-1}; q)_{n}(q^{n+1},abq^n,acq^n,adq^n,bcq^n,bdq^n,cdq^n,abcd; q)_{\infty}}$}\,\,\, \\
 &  \\
 \hline
 &  \\
 $\beta_n $ &  \begin{tabular}{c}
$\dst\frac 12\Big[-\!\frac{(1\!-\!abq^n)(1\!-\!acq^n)(1\!-\!adq^n)
(1\!-\!abcdq^{n-1})}{a(1-abcdq^{2n-1})(1-abcdq^{2n})}
\!-\!\frac{a(1\!-\!q^n)(1\!-\!bcq^{n-1})(1\!-\!bdq^{n-1})(1\!-\!cdq^{n-1})}{(1-abcdq^{2n-2})
(1-abcdq^{2n-1})}$\\[3mm] $+a\!+\!a^{-1} \Big]$
\end{tabular}\\
 &  \\
\hline
 &  \\
$\gamma_n$ & $ \dst\frac 14 \frac{(1\!-\!q^n)(1\!-\!abq^{n-1})(1\!-\!acq^{n-1})
(1\!-\!adq^{n-1})(1\!-\!bcq^{n-1})(1\!-\!bdq^{n-1})(1\!-\!cdq^{n-1})(1\!-\!abcdq^{n-2})}
{(1-abcdq^{2n-3})(1-abcdq^{2n-2})^2(1-abcdq^{2n-1})} $ \\
 &  \\
\hline
&  \\
$\,\overline{\alpha}_n\!\!=\!\widehat{\alpha}_n$ &  $4q^{-n+1}(q^{1/2}-q^{-1/2})(1-abcdq^{2n-1})$  \\&  \\
\hline &  \\
$\overline{\beta}_n(s) $ & $\begin{tabular}{r}
$-\frac{2q^{-\frac{3n}{2}+1}(q^{1/2}-q^{-1/2})(1-abcdq^{n-1})}{ab+ac+ad+bc+bd+cd}\Big[
2q^n(ab\!+\!ac\!+\!ad\!+\!bc\!+\!bd\!+\!cd)x(s\!+\!\frac n2)
\!-\!q^{n/2}(a\!+\!b\!+\!c\!+\!d)\Big]$ \end{tabular}$  \\&  \\
\hline &  \\
$\widehat{\beta}_n(s)$ & $\begin{tabular}{r}
$\overline{\beta}_n(s)-4q^{-n+1}(1-q^n)(1-abcdq^{n-1})\Delta x(s-\frac12)$ \end{tabular}$  \\&  \\
\hline
\end{tabular}  }
\end{scriptsize}
\end{center}
\end{table}
Using the identity \cite[page 156]{ran} (see also \cite[page 201]{RR})
\begin{small}
$$
(aq^{s}; q)_k(aq^{-s}; q)_k\!=\!(-1)^k a^ k q^{k(\frac{k-1}{2})}\prod_{i=0}^{k-1}
\Big[ {2 x(s)}-(aq^{i}+a^{-1}q^{-i})\Big],
$$
\end{small}%
we can rewrite \refe{pol-wil-nu} as
\[
\begin{split}\dst
P_n(x(s))_q &= \dst
\frac{(ab,ac,ad;q)_n}{(2a)^n(abcdq^{n-1};q)_n}
\sum_{k=0}^n \frac{(q^{-n},abcdq^{n-1};q)_k}{(ab,ac,ad,q;q)_k} q^ k \\
 & \qquad \times (-1)^k a^ k q^{k(\frac{k-1}{2})}
\prod_{i=0}^{k-1}\Big[ {2 x(s)}-(aq^{i}+a^{-1}q^{-i})\Big].
\end{split}
\] 
Notice that for $x(s_0)=-1$ and $x(s_1)=1$, we obtain, respectively,
\begin{equation*}
\begin{split}
P_{n}(-1)_q & = \frac{(ab,ac,ad;q)_n}{(2a)^n(abcdq^{n-1};q)_n}
{}_{4}\varphi_3 \left(\ba{c} q^{-n},abcdq^{n-1}, -a,-a  \\ ab,ac,ad  \ea
\,\bigg|\, q \,,\, q \right),\\
P_{n}(1)_q & =\frac{(ab,ac,ad;q)_n}{(2a)^n(abcdq^{n-1};q)_n}
{}_{4}\varphi_3 \left(\ba{c} q^{-n},abcdq^{n-1}, a, a \\ ab,ac,ad  \ea
\,\bigg|\, q \,,\, q \right).
\end{split}
\end{equation*}
In a similar fashion we get 
\[
\begin{split}
\Delta P_n(-1)_q=P_n(x(s_0+1))-P_n(x(s_0))&=a(1-q^{-1})(1-a)\frac{(ab,ac,ad;q)_n}{(2a)^n(abcdq^{n-1};q)_n}\\
&\times{}_{4}\varphi_3 \left(\ba{c} q^{-n},abcdq^{n-1}, aq,aq  \\ ab,ac,ad  \ea
\,\bigg|\, q \,,\, q \right),
\end{split}
\]
\[
\begin{split}
\Delta P_n(1)_q=P_n(x(s_1+1))-P_n(x(s_1))&=-a(1-q^{-1})(1+a)\frac{(ab,ac,ad;q)_n}{(2a)^n(abcdq^{n-1};q)_n}\\
&\times{}_{4}\varphi_3 \left(\ba{c} q^{-n},abcdq^{n-1}, -aq,-aq  \\ ab,ac,ad  \ea
\,\bigg|\, q \,,\, q \right).
\end{split}
\]

By inserting the values of Askey-Wilson polynomials 
given in Table \ref{tabla1} into \refe{rez12a} we arrive at 
the following identity
\bq\label{ThetaXi}
P_{n-1}(x(s))_q=\Theta(s,n)P_n(x(s))_q+\Xi(s,n) P_n(x(s+1))_q,
\eq
where
\begin{footnotesize}
\begin{equation*}
\begin{split}
\Xi(s,n)&=\dst-\frac{q^{n-1}(1-abcdq^{2n-3})(1-abcdq^{2n-2})^2}{
(1-abq^{n-1})(1-acq^{n-1})(1-adq^{n-1})(1-bcq^{n-1})(1-bdq^{n-1})(1-cdq^{n-1})}\\
&\times\frac{\Phi(s)}
{(1-abcdq^{n-2})(q^{1/2}-q^{-1/2})(1-q^n)\Delta x(s)},
\end{split}
\end{equation*}
\begin{equation*}
\begin{split}
\Theta(s,n)&=\frac{q^{n-1}(1-abcdq^{2n-3})(1-abcdq^{2n-2})^2}{
(1-abq^{n-1})(1-acq^{n-1})(1-adq^{n-1})(1-bcq^{n-1})(1-bdq^{n-1})(1-cdq^{n-1})}\\
&\times\frac{1}{(1-abcdq^{n-2})(q^{1/2}-q^{-1/2})(1-q^n)}\Bigg\{
\frac{\Phi(s)}
{\Delta x(s)}
+2q^{-n+1}(q^{1/2}-q^{-1/2})(1-abcdq^{2n-1})\\
&\times\Big[2x(s)-a-a^{-1}+(1-abcdq^{n-1})
\frac{(1-abq^n)(1-acq^n)(1-adq^n)}{a(1-abcdq^{2n-1})(1-abcdq^{2n})}+a(1-q^n)\\
&\times\frac{(1-bcq^{n-1})(1-bdq^{n-1})(1-cdq^{n-1})}{(1-abcdq^{2n-2})(1-abcdq^{2n-1})}\Big]
-\frac{2q^{-\frac{3n}{2}+1}(q^{1/2}-q^{-1/2})(1-abcdq^{n-1})}{ab+ac+ad+bc+bd+cd}\\
&\times\Big\{
2q^n(ab\!+\!ac\!+\!ad\!+\!bc\!+\!bd\!+\!cd)x(s+\mbox{$\frac n2$})
-q^{n/2}(a\!+\!b\!+\!c\!+\!d)\Big\}\\
&-4q^{-n+1}(1-q^n)(1-abcdq^{n-1})\Delta x(s-\mbox{$\frac 12$})\Bigg\}.
\end{split}
\end{equation*}
\end{footnotesize}

\section{The generalized Askey-Wilson polynomials}
In this section we consider the modification of the Askey-Wilson polynomials 
\refe{pol-wil-nu} by adding two mass points, i.e., the polynomials
orthogonal with respect to the functional $\newu=u+A\delta(x(s)-x(s_0))+B\delta(x(s)-x(s_1))$,
where $\u$ is defined in \refe{def-fun}, $x_0:=x(s_0)=-1$ and $x_1:=x(s_1)=1$.

By using \cite[\S3]{RR} the representation of the modified Askey-Wilson polynomials 
can be constructed  
\begin{equation}\label{rez1}
\vp^{A, B}(x(s))_q=P_n(x(s))_q-A\vp^{A, B}(-1)_q\K_{n-1}(x(s), -1)-B\vp^{A, B}(1)_q\K_{n-1}(x(s), 1),
\end{equation}
then the system of two equations in the two unknowns
$\vp^{A, B}(-1)_q$ and $\vp^{A, B}(1)_q$ becomes
\[\begin{split}
\vp^{A, B}(-1)_q& =P_n(-1)_q-A\vp^{A, B}(-1)_q\K_{n-1}(-1, -1))-B\vp^{A, B}(1)_q \K_{n-1}(-1, 1),\\
\vp^{A, B}(1)_q& =P_n(1)_q-A\vp^{A, B}(-1)_q\K_{n-1}(1,-1)-B\vp^{A, B}(1)_q \K_{n-1}(1, 1),
\end{split}\]
whose solution is
$$
\left(\begin{matrix} \vp^{A, B}(-1)_q \\  \vp^{A, B}(1)_q\end{matrix}\right)=
\left(\begin{matrix} 1+A\K_{n-1}(-1,-1) & B\K_{n-1}(-1,1) \\ 
A\K_{n-1}(1,-1) & 1+B\K_{n-1}(1,1)
\end{matrix}\right)^{-1}
\left(\begin{matrix} P_n(-1)_q \\ P_n(1)_q\end{matrix}\right).
$$
Notice that $\forall A,B>0$,
\bq\label{kappa(s_0,s_1)}
\begin{split}
\kappa_{n-1}(-1,1) &:= \det\left|
\begin{matrix} 1+A\K_{n-1}(-1,-1) & B\K_{n-1}(-1,1) \\ A\K_{n-1}(1,-1) & 
1+B\K_{n-1}(1,1)\end{matrix}
\right|>0.
\end{split}
\eq
Thus, by \cite[Proposition 1]{RR}
the polynomials $\vp^{A, B}(s)_q$ are well defined for all values $A,B>0$. Furthermore, 
\begin{equation}\label{rez2}
\begin{split}
\widetilde{P}_n^{A,B}(-1)_q & =\frac{(1+B\K_{n-1}(1,1))P_n(-1)_q-B
\K_{n-1}(-1,1)P_n(1)_q}{\kappa_{n-1}(-1,1)},\\
\widetilde{P}_n^{A,B}(1)_q & = \frac{(1+A\K_{n-1}(-1,-1))P_n(1)_q-A
\K_{n-1}(1,-1)P_n(-1)_q}{\kappa_{n-1}(-1,1)},
\end{split}
\end{equation}
where $\kappa_{n-1}(-1,1)$ is given in \refe{kappa(s_0,s_1)}. 

The modified Askey-Wilson polynomials satisfy the following orthogonality
relation

\[
\begin{split}
\int_{-1 }^{1} \widetilde{P}_n^{ A,B}(x)_q \widetilde{P}_m^{A,B}(x)_q &
\rho(x)dx  + A \widetilde{P}_n^{A,B}(-1)_q
\widetilde{P}_m^{A,B}(-1)_q \\ & +B\widetilde{P}_n^{A,B}(1)_q
\widetilde{P}_m^{A,B}(1)_q = \delta_{n,m} \widetilde{d}_n^2,
\end{split}
\]
where $\rho$ and $d_n$ denote 
the weight function and the norm of the Askey-Wilson polynomials (see Table
\ref{tabla1}\footnote{We have chosen $\rho(s)$ in such a way that
$\int_{x = -1 }^{1} \rho(x)dx=1$, i.e., to be a probability measure.}) and
$$
\widetilde{d}_n^2= \pe{\newu}{\p_n^2(x)}=d_n^2+ A\p_n^{A, B}(-1)_q P_n(-1)_q + B \p_n^{A, B}(1)_q P_n(1)_q.
$$

\subsection*{Representation formulas for the generalized Askey-Wilson polynomials}

Consider the representation formula \refe{rez1} where the $n$-th kernel
can be computed by the formulas \refe{rez14} and \refe{rez15}. In fact,
by using the  main datas of Askey-Wilson polynomials (see Table \ref{tabla1})
in \refe{rez15}, we obtain

\begin{equation}\label{ker-a}
\K_{n-1}(x(s), -1)=\varkappa_{-1}(s, n)
P_{n-1}(x(s))_q+\conju{\varkappa}_{-1}(s, n)
\frac{\Delta P_{n-1}(x(s))_q}{\Delta x(s)},
\end{equation}
where
\begin{equation*}
\begin{split}
 \varkappa_{-1}(s,n) \!\!& =\!\!\dst  \frac{(abcdq^{n-2};q)_{n-1}(abcd,\,q^n,\!abq^{n-1},\!acq^{n-1},\!adq^{n-1},
 \! bcq^{n-1},\! bdq^{n-1},\! cdq^{n-1};q)_{\infty}}{2^{-2n+2}q(q^{1/2}-q^{-1/2})(abcdq^{2n-3};q)_{\infty}
 (q,ab,ac,ad,bc,bd,cd;q)_{\infty}}\\
 &\times\Bigg\{\Big[\frac{q^{\frac{n+1}{2}}(q^{1/2}-q^{-1/2})(1-abcdq^{n-2})[\frac{q^{\frac{n-1}{2}}+q^{-\frac{n-1}{2}}}{2}+x(s+\frac{n-1}{2})]}
 {x(s)+1}\\
 &+\frac{(1-q^n)(1-abcdq^{n-1})\Delta x(s-\frac12)}
 {x(s)+1}\Big]P_{n-1}(-1)\\
 &+\frac{2(1+a)(1+b)(1+c)(1+d)}{(q+q^{-1}-2)[x(s)+1]}
\Delta P_{n-1}(-1)\Bigg\},
\end{split}
\end{equation*}

\begin{equation*}
\begin{split}
\conju{\varkappa}_{-1}(s,n) \!\!&=\!\!\dst \frac{(abcdq^{n-2};q)_{n-1}(abcd,\,q^n,\! abq^{n-1},\! acq^{n-1},\! adq^{n-1},
 \! bcq^{n-1},\! bdq^{n-1},\! cdq^{n-1};q)_{\infty}}{2^{-2n+4}q^{-n+2}(q^{1/2}-q^{-1/2})(abcdq^{2n-3};q)_{\infty}
 (q,ab,ac,ad,bc,bd,cd;q)_{\infty}}\\
 &\times \frac{\Phi(s)}{x(s)+1}P_{n-1}(-1),
\end{split}
\end{equation*}

\begin{equation}\label{ker-b}
\K_{n-1}(x(s), 1 )=\varkappa_{1}(s, n)P_{n-1}(x(s))_q+
\conju{\varkappa}_{1}(s, n)
\frac{\Delta P_{n-1}(x(s))_q}{\Delta x(s)},
\end{equation}
where
\begin{equation*}
\begin{split}
\varkappa_{1}(s, n)\!\!&=\!\!\dst  \frac{(abcdq^{n-2};q)_{n-1}(abcd,\,q^n,\! abq^{n-1},\! acq^{n-1},\! adq^{n-1},
 \! bcq^{n-1},\! bdq^{n-1},\! cdq^{n-1};q)_{\infty}}{2^{-2n+2}q(q^{1/2}-q^{-1/2})(abcdq^{2n-3};q)_{\infty}
 (q,ab,ac,ad,bc,bd,cd;q)_{\infty}}\\ &\times\Bigg\{\Big[-\frac{q^{\frac{n+1}{2}}(q^{1/2}-q^{-1/2})(1-abcdq^{n-2})[\frac{q^{\frac{n-1}{2}}+q^{-\frac{n-1}{2}}}{2}-x(s+\frac{n-1}{2})]}
 {x(s)-1}\\
 &+\frac{(1-q^n)(1-abcdq^{n-1})\Delta x(s-\frac12)}
 {x(s)-1}\Big]P_{n-1}(1)\\
 &-\frac{2(1-a)(1-b)(1-c)(1-d)}{(q+q^{-1}-2)[x(s)-1]}\Delta P_{n-1}(1)\Bigg\},
\end{split}
\end{equation*}
\begin{equation*}
\begin{split}
\conju{\varkappa}_{1}(s, n)\!\!&=\!\! \dst \frac{(abcdq^{n-2};q)_{n-1}(abcd,\,q^n,\! abq^{n-1},\! acq^{n-1},\! adq^{n-1},
 \! bcq^{n-1},\! bdq^{n-1},\! cdq^{n-1};q)_{\infty}}{2^{-2n+4}q^{-n+2}(q^{1/2}-q^{-1/2})(abcdq^{2n-3};q)_{\infty}
 (q,ab,ac,ad,bc,bd,cd;q)_{\infty}}\\
 &\times \frac{\Phi(s)}{x(s)-1}P_{n-1}(1).
\end{split}
\end{equation*}

By substituting \refe{ker-a} and \refe{ker-b} into \refe{rez1}, one finds
\begin{equation}\label{rez63a}
\begin{split}
\dst\widetilde{P}_n^{A, B}(x(s))_q  = \dst  P_n(x(s))_q+
\overline{A}(s, n)P_{n-1}(x(s))_q & +\overline{B}(s, n)
\frac{\Delta P_{n-1}(x(s))_q}{\Delta x(s)},
\end{split}
\end{equation}
 \begin{equation*}
\begin{split}
\overline{A}(s, n)=& -A\widetilde{P}_n^{A, B}(-1)_q\varkappa_{-1}(s, n)-
B \widetilde{P}_n^{A, B}(1)_q\varkappa_{1}(s, n),\\[4mm]
\overline{B}(s, n)=& -A\widetilde{P}_n^{A, B}(-1)_q\conju{\varkappa}_{-1}(s, n)
-B\widetilde{P}_n^{A, B}(1)_q\conju{\varkappa}_{1}(s, n),
\end{split}
\end{equation*}
where $\widetilde{P}_n^{A, B}(-1)_q$ and $\widetilde{P}_n^{A, B}(1)_q$ are given
in \refe{rez2}. 
Notice that the involved
functions $\overline{A}$ and $\overline{B}$ as well as 
${\Delta P_{n-1}(s)_q}/{\Delta x(s)}$ in \refe{rez63a} 
are not, in general, polynomials in $x(s)$. Thus,
it is not easy to see that $\widetilde{P}_n^{A, B}(s)_q$ in \refe{rez63a}
is a polynomial of degree $n$ in $x(s)$ which is even a simple consequence of \refe{rez1}. 
Notice that if we use \refe{rez14} instead of \refe{rez15} 
we can obtain a formula similar to \refe{rez63a} but in terms of the backward difference
operator.

From the Christoffel Darboux formula and the TTRR for the Askey-Wilson polynomials 
another representation formula for the modified Askey-Wilson polynomials follows
(see e.g. \S 3 in \cite{RR})
\bq \label{repfor-n-rac}
\phi(s)\widetilde{P}_n^{A, B}(x(s))_q=A(s;n)P_n(x(s))_q+B(s;n)P_{n-1}(x(s))_q,
\eq
with the coefficients
 \begin{equation}\label{ABsn}
\begin{split}
\phi(s) & =[x(s)^2-1], \\[4mm]
A(s, n)& =\phi(s) - \dst\frac{1}{d_{n-1}^2 } \Big\{
A \widetilde{P}_n^{A, B}(-1)_q P_{n-1}(-1)_q [x(s)-1] \\
&  \qquad\qquad\qquad\qquad\qquad\qquad 
+ B \widetilde{P}_n^{A, B}(1)_q P_{n-1}(1)_q[x(s)+1] \Big\}, \\
B(s, n)& = \dst\frac{1}{d_{n-1}^2 } \Big\{
A \widetilde{P}_n^{A, B}(-1)_q P_{n}(-1)_q [x(s)\!-\!1] 
+B \widetilde{P}_n^{A, B}(1)_q P_{n}(1)_q[x(s)\!+\!1] \Big\},  \\
\end{split}
\end{equation}
where $\widetilde{P}_n^{A, B}(-1)_q$ and $\widetilde{P}_n^{A, B}(1)_q$ are defined
in \refe{rez2}. 

Furthermore, there is one more representation formula
for the modified Askey-Wilson families which can be obtained by
substituting the relation \refe{ThetaXi} in \refe{repfor-n-rac}
\bq\label{repfor-n-rac2}
\phi(s)\widetilde{P}_n^{A, B}(x(s))_q=a(s;n)P_n(x(s))_q+b(s;n)P_n(x(s+1))_q,
\eq
where $a(s;n)=A(s;n)+B(s;n)\Theta(s;n)$, $b(s;n)=B(s;n)\Xi(s;n)$,
and $A$, $B$ and $\Theta$, $\Xi$ are given by \refe{ABsn} and \refe{ThetaXi}, respectively.

If we change in \refe{repfor-n-rac2} $s$ by $s+1$ and  $s$ by $s-1$ and then use \refe{rez8} to eliminate
$P_n(x(s+2))_q$ and $P_n(x(s-2))_q$, respectively, we obtain 
\begin{equation}\label{5} u(s)\vp(x(s+1))_q=c(s,n)P_n(x(s))_q+d(s,n)P_n(x(s+1))_q,
\end{equation}
and
\begin{equation}\label{7}
v(s)\vp(x(s-1))_q=e(s,n)P_n(x(s))_q+f(s,n)P_n(x(s+1))_q,
\end{equation} 
respectively, where $u(s)=A_{s+1}\phi(s+1)$, $c(s,n)=-C_{s+1}b(s+1,n)$, and
$d(s,n)=A_{s+1}a(s+1,n)-b(s+1,n)(\lambda_n+B_{s+1})$,  
$v(s)= C_s\phi(s-1)$, $e(s,n)=C_s b(s-1,n)-a(s-1,n)(\lambda_n+B_s)$, and
$f(s,n)=-A_s a(s-1,n)$.
Then (\ref{repfor-n-rac2}), (\ref{5}) and (\ref{7}) lead to
\begin{equation}\label{9}
\left|\begin{matrix}
\phi(s)\vp(x(s))_q & a(s,n) & b(s,n) \\
u(s)\vp(x(s+1))_q & c(s,n) & d(s,n) \\
v(s)\vp(x(s-1))_q & e(s,n) & f(s,n)
\end{matrix}\right|=0.
\end{equation}
Expanding the determinant \refe{9} by the first column, we get the following
second order linear difference equation for $\widetilde{P}_n^{A, B}(x(s))_q$
\begin{equation}\label{10}
\widetilde{\phi}(s,n)\vp(x(s-1))_q+\widetilde{\varphi}(s,n)\vp(x(s))_q+\widetilde{\xi}(s,n)\vp(x(s+1))_q=0,\\
\end{equation}
\begin{equation*}
 \begin{split}
\widetilde{\phi}(s, n)& =v(s)\Big[a(s, n)d(s, n)-b(s, n)c(s, n)\Big], \\
\widetilde{\varphi}(s, n)&=\phi(s)\Big[c(s, n)f(s, n) - d(s, n)e(s, n)\Big],\\
\widetilde{\xi}(s, n)&=u(s)\Big[b(s, n)e(s, n)-a(s, n)f(s, n)\Big].
\end{split}
\end{equation*}
Thus, the generalized Askey-Wilson polynomials satisfy 
a second order linear difference equation \refe{10}
with polynomial coefficients which explicitly depend on $n$.

Moreover one can obtain the TTRR of the monic generalized
Askey-Wilson polynomials with two mass points (for details see Eqs. (20), (21) in \cite{RR}) 
$$
x(s)\widetilde{P}_n^{A, B}(x(s))_q=\widetilde{P}_{n+1}^{A, B}(x(s))_q
+\widetilde{\beta}_n\widetilde{P}_n^{A, B}(x(s))_q+\widetilde{\gamma}_n\widetilde{P}_{n-1}^{A, B}(x(s))_q, 
\quad n\in\mathbb{N},
$$
where
\begin{small}
\begin{equation*}
\begin{split}
\widetilde{\beta}_n &=   \beta_n
-A \left(\frac{\widetilde{P}_n^{A,B}(-1)_q P_{n-1}(-1)_q}{{d}_{n-1}^2}-
\frac{\widetilde{P}_{n+1}^{A,B}(-1)_qP_{n}(-1)_q }{{d}_{n}^2}\right) \\
&-B \left(\frac{\widetilde{P}_n^{A,B}(1)_q P_{n-1}(1)_q }{{d}_{n-1}^2}-
\frac{\widetilde{P}_{n+1}^{A,B}(1)_q P_{n}(1)_q }{{d}_{n}^2}\right), \\
\widetilde{\gamma}_n &= \gamma_n\frac{1\!+\!\Delta_n^{A,B}}{1\!+\!\Delta_{n-1}^{A,B}},
\,\, \Delta_n^{A,B}=
\frac{ A \widetilde{P}_n^{A,B}(-1)_q P_{n}(-1)_q }{{d}_{n}^2}\!+\!
\frac{B \widetilde{P}_n^{A,B}(1)_q P_{n}(1)_q}{{d}_{n}^2}.
\end{split}
\end{equation*}
\end{small}


\subsubsection*{Representation of $\widetilde{P}_n^{A,B}(x(s))_q$ in  terms of
basic series}
In this section, we obtain an explicit formula for 
$\widetilde{P}_n^{A,B}(x(s), a, b,c,d|q)$
in terms of basic hypergeometric series. In fact,
substituting \refe{pol-wil-nu} into \refe{repfor-n-rac} we obtain
\[
\ba{rl}
\phi(s)\widetilde{P}_n^{A, B}(x(s))_q \!\!=\!\!\! &\!\!
\dst\frac{(ab,ac,ad;q)_{n-1}}{(2a)^{n-1}(abcdq^{n-2};q)_{n-1}}
\dst\!\! \sum_{k=0}^{\infty}\frac{(q^{-n}, abcdq^{n-2}, aq^{s}, aq^{-s}; q)_k}
{(ab,ac, ad, q; q)_k}q^k
\Pi_1(q^k),
\ea
\]
where $\phi(s)$, $A(s, n)$ and $B(s, n)$ are given in \refe{ABsn} and
\bq\label{Pi_1}
\begin{split}
\Pi_1(q^k)&  = 
A(s, n)\frac{(1-abq^{n-1})(1-acq^{n-1})(1-adq^{n-1})(1-abcdq^{n+k-2})}
{2a(1-abcdq^{2n-3})(1-abcdq^{2n-2})}\\
& \quad + B(s, n)\frac{(1-q^{-n+k})}{(1-q^{-n})}\\
& = -\Big\{A(s, n) abcdq^{n-2}
\vartheta^{a,b,c,d}_n+B(s, n)q^{-n}\Big\}\frac{(q^k-q^{\kappa(s)})}{1-q^{-n}},
\end{split}
\eq
where
\[\begin{split}
q^{\kappa(s)} & =\frac{A(s, n) \vartheta^{a,b,c,d}_n+ B(s, n)}
{A(s, n) abcdq^{n-2}\vartheta^{a,b,c,d}_n+ B(s, n)q^{-n}}, \\
\vartheta^{a,b,c,d}_n & =\frac{(1-abq^{n-1})(1-acq^{n-1})(1-adq^{n-1})(1-q^{-n})}
{2a(1-abcdq^{2n-3})(1-abcdq^{2n-4})}.
\end{split}
\]
By taking into account the identity $(q^k-q^{z})(q^{-z};q)_k=
(1-q^{z})(q^{1-z}; q)_k$ we obtain
\bq\label{hypreprac2}
\begin{split}\dst
\phi(s)\widetilde{P}_n^{A, B}(x(s))_q =& {D}_n^{a,b,c,d}(s)
{}_{5}\varphi_4 \! \left(\!\ba{c} q^{-n},abcdq^{n-2}, aq^{s},
aq^{-s}, q^{1-\kappa(s)} \\ ab,ac,ad, q^{-\kappa(s)}\ea\!\!
\bigg|q, q \!\right),
\end{split}
\eq
where
\begin{equation*}
\begin{array}{rl}
{D}_n^{a,b,c,d}(s)\!\!&\!=\!\!\dst\dst\frac{-(ab,ac,ad;q)_{n-1}}{(2a)^{n-1}(abcdq^{n-2};q)_{n-1}}
\frac{1\!-\!q^{\kappa(s)}}{1\!-\!q^{-n}}\Big\{A(s, n) abcdq^{n-2}\vartheta^{a,b,c,d}_n\!+\!
B(s, n)q^{-n}\Big\}.
\end{array}
\end{equation*}

\begin{rem} Notice that $\phi(s)\widetilde{P}_n^{A, B}(x(s))_q$, in the left hand side of \refe{hypreprac2}, 
is a polynomial of degree $n+2$ in $x(s)$ which follows from \refe{repfor-n-rac}
and \refe{ABsn}. In order to see that formula \refe{hypreprac2} gives a polynomial of degree $n+2$ it is 
sufficient to notice that the function $\Pi_1$ defined in \refe{Pi_1} is a polynomial in $x(s)$,
which follows from that fact that $A(s,n)$ and $B(s,n)$ are polynomial of degree 
2 and 1 in $x(s)$, respectively, \refe{ABsn}. 
\end{rem}   



\begin{rem}
We note that the properties of the modified Askey-Wilson polynomials
with one mass point at $x=\pm1$ can be obtained from
the ones with two mass points by putting $A=0$ or $B=0$, 
respectively.
\end{rem}

\begin{rem}
We remark that the relation between Askey-Wilson $P_n(x;a,b,c,d;q)$ polynomials 
 defined in \refe{pol-wil-nu}
and $q$-Racah polynomials $u_n^{\alpha,\beta}(\mu(t),\widetilde{a},\widetilde{b})$ defined on the lattice
$\mu(t)=[t]_q[t+1]_q=c_1(q^t+q^{-t-1})+c_3$, $c_1=q^{1/2}(q^{1/2}-q^{-1/2})^{-2}$
and $c_3=-q^{-1/2}(1+q)(q^{1/2}-q^{-1/2})^{-2}$ follows 
\[
\begin{split}
\frac{2^n}{(q^{1/2}-q^{-1/2})^{2n}}P_n(2c_1q^{-1/2}x+c_3,q^{\widetilde{a}+\frac 12},
q^{\beta-\widetilde{a}+\frac 12},q^{\alpha+\widetilde{b}+\frac 12},q^{-\widetilde{b}+\frac 12};q)=
u_n^{\alpha,\beta}(\mu(t),\widetilde{a},\widetilde{b})
\end{split}
\]
by setting $e^{2i\theta}=q^{2s}=q^{2t+1}$, where
\[
\begin{split}
u_n^{\alpha,\beta}(\mu(t),\widetilde{a},\widetilde{b})&=q^{-\frac n2(2\widetilde{a}+1)}
\frac{(q^{\widetilde{a}-\widetilde{b}+1},q^{\beta+1},q^{\widetilde{a}+\widetilde{b}+\alpha+1};q)_n}
{(q^{1/2}-q^{-1/2})^{2n}(q^{\alpha+\beta+n+1};q)_n}\\
&\times{}_{4}\varphi_3 \left(\ba{c} q^{-n},q^{\alpha+\beta+n+1}, q^{\widetilde{a}-t},
q^{t+\widetilde{a}+1}  \\ q^{\widetilde{a}-\widetilde{b}+1},q^{\beta+1},q^{\widetilde{a}+\widetilde{b}+\alpha+1}  \ea
\,\bigg|\, q \,,\, q \right).
\end{split}
\]
\end{rem}
Combining the above limit, with the ones  considered in
\cite[\S 4.2]{RR} we can construct the analog of $q$-Askey
Tableau for the Krall-type polynomials. For more details 
on how one should take the limits we refer to the paper 
\cite{stk98}.

\section*{Concluding remarks} 
In this paper we have constructed a generalized Askey-Wilson 
polynomials by adding two mass points at the end of the interval of 
orthogonality and obtained some of their properties, as the 
TTRR and the representation as basic hypergeometric series.
In particular, we have showed that they satisfy a second order linear
$q$-difference equation on the lattice $x(s)=(q^s+q^{-s})/2$
(see \refe{10}). This equation has the form \refe{rez8}
but with coefficients that explicitely depend on $n$, the 
degree of the polynomials. In general, they will not satisfy 
a higher order difference equation with coefficients independent 
of $n$. An example of such polynomials satisfying
a higher order difference equation with coefficients independent 
of $n$ was constructed in \cite{VZ01}.

\section*{Acknowledgements:}  We want to thank the
unknown referees for their suggestions that helped us to improve the paper
and for pointing out the paper \cite{VZ01}. 
This work was partially supported by MTM2009-12740-C03-02
(Ministerio de Econom\'\i a y Competitividad), FQM-262, FQM-4643, FQM-7276 (Junta de Andaluc\'\i a),
Feder Funds (European Union).
The second author is also supported by a grant from T\"{U}B\.{I}TAK, the Scientific and
Technological Research Council of Turkey. She also thanks to the Departamento 
de An\'alisis Matem\'atico and IMUS for their kind hospitality.



\end{document}